\documentclass{article}
\usepackage{amsmath, amsthm}

\newcommand{\mpt}{\;\text{mod}\,p^3}
\newcommand{\mpd}{\;\text{mod}\,p^2}
\newcommand{\mpu}{\;\text{mod}\,p}
\newcommand{\mb}{\mathcal{B}}
\newcommand{\se}{\sum_{\begin{array}{l}k=0\\\text{$k$ even}\end{array}}^{p-1}H_k}
\newcommand{\so}{\sum_{\begin{array}{l}k=0\\\text{$k$ odd}\end{array}}^{p-1}H_k}
\newcommand{\hp}{H^{'}_{p-1}}
\newcommand{\wh}{\widehat}
\newcommand{\bc}{$\Large{\times}$}
\newcommand{\mg}{\mathcal{G}}
\newcommand{\wih}{\widehat}

\author{Claire Levaillant}
\title{Congruences related to Eulerian and Euler~numbers}
\begin{document}
\maketitle
\begin{center}\textit{Dedicated to my late grandmother M-F. Andr\'e}\end{center}
\begin{abstract} Let $H_k$ denote the harmonic number of order $k$ and $\hp$ denote the sum of order $p-1$ of the odd reciprocals of integers. Given an odd prime number $p>3$, by using in particular some congruences on Eulerian numbers, we show the following congruences:
\begin{equation*}\sum_{k=0}^{p-1}(-1)^k\,kH_k=\frac{1}{2}\,H^{'}_{p-1}\;\mpu\qquad
\sum_{\begin{array}{l}k=0\\\text{$k$ even}\end{array}}^{p-1}kH_k=\frac{1}{4}(H^{'}_{p-1}-1)\mpu\end{equation*}
\begin{equation*}\sum_{k=1}^{p-1}\frac{1}{k2^k}=\hp\mpu\end{equation*}
\indent Let $T_k$ denote the $k$-th tangent number defined by the series expansion $\tan\,z=\sum_{k\geq 0}T_k\frac{z^k}{k!}$ for $|z|<\frac{\pi}{2}$. Given an old prime number $p>3$, by using in particular some congruence on the Euler numbers, we show the congruence:
\begin{equation*}
\sum_{k=1}^{\frac{p-1}{2}}(-1)^k(2k-1)T_{2k-1}=\begin{cases}0\;\;\mpu&\text{if $p\equiv\,1$ mod 4}\\
2\;\;\mpu&\text{if $p\equiv\,3$ mod 4}\end{cases}
\end{equation*}
%\indent Let $G_k$ denote the $k$-th Genocchi number. Given an odd prime $p$, by using in particular the Mirimanoff congruences, we show the congruence
%\begin{equation*} \sum_{k=}^{}\frac{G_k}{k2^k}=\mpu\end{equation*}
\end{abstract}
\section{Introduction}
The goal of this paper is to derive some nice congruences from congruences on Eulerian numbers on one hand and from congruences on Euler numbers on the other hand. Some of these congruences were also obtained by other authors by other means. We focus here on these two types of numbers, namely the Eulerian numbers and the Euler numbers, and on old and more recent progress that were made on their arithmetic properties. These numbers also relate to the Bernoulli numbers \cite{NIEL} and to the Genocchi numbers \cite{GENO}. Though originally defined through coefficients in the Taylor series expansions of some classical functions, they both have a really neat combinatorial interpretation.

\subsection{An introduction to the Eulerian numbers}

We begin introducing the Eulerian numbers. These numbers first appear in Euler's book dating from $1755$ \cite{EULE} pp $487-491$. Euler investigated polynomials that are shifted forms of what are now called Eulerian polynomials. \\The generating function given by Euler is:
$$\frac{1-x}{e^t-x}=\sum_{n=0}^{\infty}H_n(x)\frac{t^n}{n!}$$
The Eulerian numbers are then defined by:
$$(x-1)^n\,H_n(x)=\sum_{m=0}^{n-1}E(n,m)x^m$$
MacMahon was first to provide in \cite{MACM} a combinatorial interpretation for the Eulerian numbers.
The Eulerian number $E(n,m)$ counts the number of permutations of $Sym(n)$ with $m$ ascents. If $(i_1,i_2,\dots,i_n)$ denotes the permutation mapping the integer $j$ onto $i_j$, an ascent (resp a descent) is when $i_{k+1}>i_k$ (resp $i_{k+1}<i_k$). The Eulerian numbers can be calculated using a recursive formula:
$$E(n,m)=(n-m)E(n-1,m-1)+(m+1)E(n-1,m)\qquad (R_{n,m})$$
A concrete interpretation of the Eulerian numbers is a special case of Simon Newcomb's problem. The latter problem is associated with a game of patience played with a deck of cards numbered from $1$ to $n$. Turn over the top card and place it down face up. At each step, turn over the next card and place it on top of the previous card if its number is less than the number of the previous card or start a new pile in the opposite situation. Then, $E(n,m)$ is the number of arrangements of the deck leading to $m+1$ piles. This is only a special case of the Simon Newcomb's problem. In the general Simon Newcomb's problem, several cards can have the same value. When the deck has $n$ kinds of cards and $i$ cards of each kind, the answer is $E^{(i)}(n,m)$ where $E^{(i)}(n,m)$ can be defined by a generalization of the Worpitsky \cite{WORP} identity due to Shanks \cite{SHAN}:
\begin{eqnarray*}
x^n&=&\sum_{m=0}^{n-1}E(n,m)\binom{x+m}{n}\qquad\textit{Worpitsky identity \cite{WORP}}\\
\binom{x}{i}^n&=&\sum_{m=0}^{i(n-1)}E^{(i)}(n,m)\binom{x+m}{in}\;\,\,\textit{Shanks identity}\;\;\;\;\,\textit{\cite{SHAN}}
\end{eqnarray*}
Because the Eulerian numbers form a triangular array called Euler's triangle which presents the same symmetries as the Pascal triangle (this is a pretty straightforward fact from their combinatorial definition), these numbers are sometimes called the "triangular numbers".
The first few Eulerian numbers get displayed on the triangular array below.
\begin{center}\begin{tabular}{|c|c|c|c|c|c|c|c|}\hline
n\textbackslash m&0&1&2&3&4&5&6\\\hline
1&1&&&&&&\\\hline
2&1&1&&&&&\\\hline
3&1&4&1&&&&\\\hline
4&1&11&11&1&&&\\\hline
5&1&26&66&26&1&&\\\hline
6&1&57&302&302&57&1&\\\hline
7&1&120&1191&2416&1191&120&1\\\hline
\end{tabular}
\end{center}
$$\begin{array}{l}\end{array}$$

\noindent An immediate consequence readable on the table
is that an alternating sum of Eulerian numbers is zero in the case when $n$ is even. In the case when $n$ is odd, the alternating sums of Eulerian numbers form a sequence of integers with alternating signs. We will refer to this sequence again in $\S\,1.2$.\\
Another interesting feature of these alternating sums of Eulerian numbers is that they relate to the Bernoulli numbers by:
\begin{equation*}
2^{n+1}(2^{n+1}-1)\frac{B_{n+1}}{n+1}=\sum_{m=0}^{n-1}(-1)^m\,E(n,m)\qquad (F)_n
\end{equation*}
A closed formula for the Eulerian numbers appears in the English version of Louis Comtet's book \cite{COMT} on page $243$. It is the following:
$$E(n,m)=\sum_{k=0}^{m+1}(-1)^k\binom{n+1}{k}(m+1-k)^n$$
When $p$ is an odd prime, some congruences modulo powers of $p$ on the Eulerian numbers of the form $E(p-2,m)$ can be easily derived from this closed formula and the results of $\S\,2$ are in part based on these congruences.\\
Other types of congruences on the Eulerian numbers appear in $1952$ in the works of Carlitz and Riordan, see $\S\,3$ of \cite{CARI}. They show some periodicity properties: for fixed $m$, if $i=\lceil log_p(m+1)\rceil$ and $j\leq n$, then
\begin{eqnarray*}E(n,m)&=&E(n+p^{i+j-1}(p-1),m)\;\text{mod}\,p^j\\
E^{(i)}(n,m)&=&E(n+p^{i+j-1}(p-1),m)\;\text{mod}\,p^j\end{eqnarray*}
Eulerian numbers appear in a variety of contexts, including when dealing with convolutions. In \cite{GOUL}, the author studies sums of convolved powers and shows that:
$$\sum_{k=0}^{n}k^i(n-k)^j=\sum_{r=0}^{i+j}\binom{i+j+1+n-r}{i+j+1}\sum_{s=0}^r E(i,s)E(j,r-s)$$
He also comes up with a formula involving Stirling numbers of the second kind.
In \cite{LEV2}, we use some congruences on the Eulerian numbers $E(p-2,m)$ in order to derive a congruence on a convolution of order $(p-1)$ of powers of two weighted Bernoulli numbers with Bernoulli numbers in terms of harmonic numbers and generalized harmonic numbers.

The next paragraph is concerned with introducing the Euler numbers. It follows the same sketch as $\S\,1.1$ in that we first introduce these numbers by series expansion, then deal with the combinatorial aspects. At the end of the paragraph, we present some closed forms and survey the currently known arithmetic properties.
\subsection{An introduction to the Euler numbers}
We now introduce the Euler numbers. Because different authors use different notations, we will cautiously settle the terminology here, in a way which may embrace the totality of the possibilities. The most common definition of the Euler numbers in the literature is the following. This sequence of numbers usually denoted by $(E_n)_{n\geq 0}$ is defined by the Taylor series expansion
$$\frac{1}{cosh\,t}=\frac{2e^t}{e^{2t}+1}=\sum_{n=0}^{\infty}\frac{E_n}{n!}t^n\qquad |t|<\frac{\pi}{2}$$
Consequently,
\begin{equation*}
E_0=1\qquad E_{2n-1}=0\qquad \sum_{s=0}^{n}\binom{2n}{2s}E_{2s}=0\qquad\forall\,n\geq 1,
\end{equation*}
sometimes abbreviated in the literature by $E_0=1$ and $(E+1)^n+(E-1)^n=0$, using the symbolic notation $E^n:=E_n$.\\
The first few Euler numbers get provided in the table below. \\\\

\begin{tabular}{|c|ccccccccccccc|}\hline
$n$&0&1&2&3&4&5&6&7&8&9&10&11&12\\\hline
$E_n$&1&0&-1&0&5&0&-61&0&1385&0&-50521&0&2702765\\\hline\end{tabular}
$$\begin{array}{l}\end{array}$$
%E_0=1,\qquad E_2=-1,\qquad E_4=5,\qquad E_6=-61
The Euler numbers relate to a special value of the Euler polynomials by
$$E_n=2^nE_n\left(\begin{array}{l}\frac{1}{2}\end{array}\right),$$
where $$\frac{2e^{xt}}{e^t+1}=\sum_{n=0}^{\infty}E_n(x)\frac{t^n}{n!}\qquad |t|<\pi$$
They satisfy to the recursive formula (see \cite{CHCW}):
\begin{equation*}
E_n+2^{n-1}\sum_{k=0}^{n-1}\binom{n}{k}\frac{E_k}{2^k}=1\qquad(n\geq 1)
\end{equation*}
Viennot's approach in \cite{VIEN} is combinatorial. Thus, he only deals with sequences of positive integers, see for instance his table on page $5$. His definition of the Euler numbers is in part for that purpose slightly different. He defines the Euler number $E_n$ as the tangent number, say $T_n$ when $n$ is odd (see e.g. \cite{KNBU}) and the secant number, say $S_n$ when $n$ is even. Both integers are always positive and
\begin{equation*}
sec(t)=\sum_{k=0}^{\infty}(-1)^kE_{2k}\frac{t^{2k}}{(2k)!}\qquad |t|<\frac{\pi}{2}
\end{equation*}
Therefore, $$S_{2k}=(-1)^kE_{2k}$$
The main advantage of this definition is to come up with a common combinatorial interpretation for the Euler number $E_n$ as the number of alternating permutations of $Sym(n)$, see Proposition $3.3$ on page $24$ of \cite{VIEN}. Indeed, within this combinatorial definition, we can express $2E_{n+1}$ as a binomial convolution:
$$2E_{n+1}=\sum_{k=0}^n\binom{n}{k}E_kE_{n-k}\qquad n\geq 1\qquad(C_n)$$
Setting $y=\sum_{n\geq 0}E_n\,t^n/n!$ and taking into account the initial conditions $E_0=E_1=1$, Equation $(C_n)_{n\geq 1}$ becomes the differential equation $$2y'=y^2+1,\;\; y(0)=1,$$
whose unique solution is $y=sec\,t+tan\,t$, see e.g. \cite{STAN}.\\
A permutation is "alternating" if all the descents occur exactly on all the odd integers. Alternating permutations were first studied by French combinatorialist D\'esir\'e Andr\'e in \cite{AND1} and \cite{AND2}.\\
In order to integrate Viennot's viewpoint into our own exposition, we will introduce our own generalization of the Euler numbers as originally defined by the Taylor expansion of the hyperbolic cosine. We define the generalized Euler number $\widehat{E_n}$ as:
$$\widehat{E_n}=\begin{cases}E_n&\text{if $n$ is even}\\
\sum_{m=0}^{n-1}(-1)^nE(n,m)&\text{if $n$ is odd}\end{cases}$$
Thus, the generalized Euler number $\wh{E_n}$ is the Euler number $E_n$ when $n$ is even and an alternating sum of Eulerian numbers when $n$ is odd.\\

The table below summarizes the first values for $\wh{E_n}$. \\\\

\hspace{-1.3cm}\begin{tabular}{|c|ccccccccccccc|}\hline
$n$&0&1&2&3&4&5&6&7&8&9&10&11&12\\\hline
$\wh{E_n}$&1&1&-1&-2&5&16&-61&-272&1385&7936&-50521&-353792&2702765\\\hline\end{tabular}
$$\begin{array}{l}\end{array}$$

The difference with respect to Viennot's definition is that we now take signs into account. Indeed, with this definition we have:
$$|\wh{E_n}|=\begin{cases}S_n&\text{when $n$ is even}\\
T_n&\text{when $n$ is odd}
\end{cases}$$

In order to see the latter fact, it is a good time to introduce the Genocchi numbers.
The Genocchi numbers, named after the Italian mathematician Angelo Genocchi $(1817-1889)$ is the sequence of numbers defined by the series expansion:
$$\frac{2t}{e^t+1}=\sum_{n=1}^{\infty}G_n\frac{t^n}{n!}\qquad |t|<\pi$$
They thus relate to the Bernoulli numbers by:
$$G_n=2(1-2^n)B_n$$
Suppose $n$ is odd. We see with the latter identity combined with the formula $(F)_n$ provided in $\S\,1.1$ that:
$$\wh{E_n}=-\frac{G_{n+1}}{n+1}2^n$$
Then, the $(2n+1)$-th tangent number relates to the Genocchi number $G_{2n+2}$ by
$$T_{2n+1}=|\wh{E_{2n+1}}|=\frac{|G_{2n+2}|}{n+1}2^{2n}$$
This is precisely identity $(1.6)$ on page $4$ of \cite{VIEN}. Thus, $(|\wh{E_n}|)_{n\geq 0}$ is the sequence of positive integers as defined by Viennot. \\
The first few Genocchi numbers are $1,-1,0,1,0,-3,0,17,0,-155,0,2073$.
There exists a combinatorial interpretation for the unsigned Genocchi numbers $|G_{2n}|$, provided by Dominique Dumont in \cite{DUMO}. Dumont's pioneering combinatorial interpretation is the following. In the same way $|\wh{E_n}|$ counts the number of permutations of $Sym(n)$ ascending or descending following the parity of $i$, the number $|G_{2n}|$ counts the number of permutations $\sigma$ of $Sym(2n-1)$ ascending or descending following the parity of $\sigma(i)$. Namely, $|G_{2n}|$ is the number of permutations $\sigma$ of $Sym(2n-1)$ with
$$\begin{array}{cc}
\sigma(i)<\sigma(i+1)&\text{if $\sigma(i)$ is odd}\\
\sigma(i)>\sigma(i+1)&\text{if $\sigma(i)$ is even}
\end{array}$$
Another combinatorial interpretation for the unsigned Genocchi number $|G_{2n}|$ can be found in \cite{VIEN} pp. $44-45$. Following the definition of \cite{VIEN}, a "gun" g on $[\!|1,n|\!]$ is a map from $[\!|1,n|\!]$ to $[\!|1,n|\!]$ such that $g(2i)\leq i$ and $g(2i-1)\leq i$ for each $i\geq 1$ (I personally rather view it as a pipe). A gun is said to be alternating if it descends on the odd integers and ascends strictly on the even integer. The notion is best illustrated on the following array serving as an example.
\begin{center}\begin{tabular}{|c|c|c|c|c|c|c|c|}\cline{7-8}
\multicolumn{6}{c|}{}&$\Large{\times}$&\\\cline{5-8}
\multicolumn{4}{c|}{}&\bc&\bc&&\\\cline{3-8}
\multicolumn{2}{c|}{}&\bc&&&&&\bc\\\hline
\bc&\bc&&\bc&&&&\\\hline
\end{tabular}\end{center}
$$\begin{array}{l}\end{array}$$
By Proposition $4.11$ on page $45$ of \cite{VIEN}, the unsigned Genocchi number $|G_{2n}|$ is the number of alternating guns on $[\!|1,2n-2|\!]$.

Like for the Eulerian numbers, there also exist some closed formulas for the Euler numbers with even indices. The dedicated reference on the topic is \cite{WEIQ}. Amongst their results obtained in year $2015$, the authors show that $E_{2n}$ can be expressed as the determinant of a matrix of size $2n$, up to a sign, or as a double sum, namely:
\begin{eqnarray*}E_{2n}&=&(-1)^n\Big|\binom{i}{j-1}cos\big((i-j+1)\frac{\pi}{2}\big)\Big|_{(2n)\times(2n)}\\
&=&(2n+1)\sum_{l=1}^{2n}(-1)^n\frac{1}{2^l(l+1)}\binom{2n}{l}\sum_{q=0}^l\binom{l}{q}(2q-l)^{2n}\end{eqnarray*}
Some of their formulas also involve the Stirling numbers of the second kind. \\A decade before in $2004$, while establishing congruences for the Euler numbers, Kwang-Wu Chen had also found a closed formula involving Genocchi numbers, namely (see \cite{CHEN}):
$$E_n=1+\frac{1}{n+1}\sum_{k=2}^{n+1}\binom{n+1}{k}2^{k-1}G_k\qquad (n\geq 0)$$
Like for the Eulerian numbers, the arithmetic study of the Euler numbers is rather old, but a resurgence of interest is only a matter of the past two decades, starting again with the work of Wenpeng Zhang \cite{ZHAN} who showed for $p$ odd prime that:
\begin{equation*}
E_{p-1}=\begin{cases}
0\;\mpu &\text{if $p\equiv 1\,\text{mod}\,4$}\\
2\;\mpu &\text{if $p\equiv 3\,\text{mod}\,4$}
\end{cases}
\end{equation*}
This is a result which pre-existed in George Ely's paper from $1880$, see \cite{GELY} p. $341$ and had been generalized by Nielsen to the even subscripts divisible by $p-1$, see \cite{NIEL} p. $273$, and generalized further to the positive integers $n$ such that $\varphi(p^a)|n$ with $a>0$ and $\varphi$ Euler's totient function, by Leonard Carlitz in \cite{CARL}. The latter author namely proves that:
\begin{equation*}
E_n=\begin{cases}
0\;\text{mod}\,p^{a+1}&\text{if $p\equiv 1\,\text{mod}\,4$}\\
2\;\text{mod}\,p^{a+1}&\text{if $p\equiv 3\,\text{mod}\,4$}
\end{cases}
\end{equation*}
His result can be viewed as an analogue of John Adams'theorem from $1878$ for the Bernoulli numbers which says that given an odd prime $p>3$,
if $p^l|n$ and $p-1\not|\,n$, then $p^l|B_n$ (see e.g. reference \cite{JOHN}).\\
In $2002$, using some results by Wells Johnson from \cite{JOHN}, Samuel Wagstaff re-discover this result in \cite{WAGS}.
Wagstaff encounters these congruences while tackling the problem of determining the prime divisors of the Euler numbers. The prime factors of the Euler numbers determine the structure of certain cyclotomic fields, see \cite{ERME}. \\
In $2004$, Kwang-Wu Chen proves an analogue of the famous Kummer type congruences for the Bernoulli numbers: given $a$ and $k$ some positive integers,
\begin{equation*}
E_{k\varphi(p^a)+2n}=\big(1-(-1)^{\frac{p-1}{2}}p^{2n}\big)E_{2n}\;\text{mod}\,p^a\qquad (n\geq 0)
\end{equation*}
Finally, the latest results on the topic date from $2008$ and are due to Yuan He and Qunying Liao in \cite{HELI}. The authors provide in particular an interesting lemma giving some general congruences on the Euler numbers. They show the following elegant congruences which hold for any odd integer $m\geq 1$:
\begin{equation*}
E_n=\sum_{l=0}^{m-1}(-1)^l(2l+1)^n\;\;(mod\;m)\qquad\qquad\qquad(n\geq 0)
\end{equation*}

Our $\S\,3$ is mostly based on congruences derived from the closed formula of Kwang-Wu Chen and also uses George Stetson Ely's pioneering congruences, later rediscovered by Wenpeng Zhang.

Last, we end this arithmetical discussion by mentioning some connections with the Fermat quotients.
Like is visible on the closed form for the Eulerian numbers provided in $\S\,1.1$, the Eulerian number $E(p-1,k)$ relates to the Fermat quotients $q_2(p),\dots,q_{k+1}(p)$. The Euler numbers are also connected to the Fermat quotients. Congruences involving Euler numbers, Fermat quotients and harmonic sums or generalized harmonic sums were first pointed out by Emma Lehmer in \cite{LEHM} and later appear again in the work of Zhi-Hong Sun in \cite{SUNH}. As part of his many results, Zhi-Hong Sun obtains in his respective Theorem $3.2$ and $3.7$ a series of congruences, one of which reads for instance:
\begin{equation*}
%\sum_{1\leq k<\frac{p}{4}}\frac{1}{k}&=&-3q_2(p)+p\big(\frac{3}{2}q_2(p)^2+(-1)^{\frac{p-1}{2}}(E_{2p-4}-2E_{p-3})\big)-p^2\big(q_2(p)^3+\frac{7}{12}B_{p-3}\big)\mpt\\
\sum_{\frac{p}{4}<k<\frac{p}{2}}\frac{1}{k}=q_2(p)-p\Big(\frac{1}{2}q_2(p)^2+(-1)^{\frac{p-1}{2}}(E_{2p-4}-2E_{p-3})\Big)+\frac{1}{3}p^2q_2(p)^3\mpt
%\sum_{\frac{p}{2}<k<\frac{3p}{4}}\frac{1}{k}&=&-q_2(p)+p\big(\frac{1}{2}q_2(p)^2+(-1)^{\frac{p-1}{2}}(6E_{p-3}-3E_{2p-4})\big)-\frac{1}{3}p^2(q_2(p)^3+14B_{p-3})\mpt\\
%\sum_{\frac{3p}{4}<k<p}\frac{1}{k}&=&3q_2(p)-p\big(\frac{3}{2}q_2(p)^2+(-1)^{\frac{p-1}{2}}(6E_{p-3}-3E_{2p-4})\big)+p^2\big(q_2(p)^3+\frac{59}{12}B_{p-3}\big)\mpt
\end{equation*}
The other ranges for $k$, namely $1\leq k<\frac{p}{4}$ (leading to a generalization of Emma Lehmer's original congruence modulo $p$ of \cite{LEHM}),              $\frac{p}{2}<k<\frac{3p}{4}$ and $\frac{3p}{4}<k<p$ get treated by Zhi-Hong Sun as well. In \cite{CAIF}, the authors also find some generalization of Emma Lehmer's original congruence by showing that:
\begin{equation*}
\sum_{\begin{array}{l}r=1\\p\not|r\end{array}}^{\lfloor p^l/4\rfloor}\frac{1}{r^2}\equiv (-1)^{(p^l-1)/2}\,4E_{\varphi(p^l)-2}\begin{cases}(\text{mod}\,p^l)&\text{for $p\geq 5$}\\(\text{mod}\,3^{l-1})&\text{for $p=3$}\end{cases}
\end{equation*}
Later, in \cite{COSD}, the authors extend these congruences to arbitrary moduli $n$.\\
In an unrelated work in \cite{JAKU}, by using some Galois theory, Stanislav Jakubec obtains a congruence modulo $p^3$ among the Euler numbers $E_{p-1}$ and $E_{2p-2}$ and the Fermat quotients $q_2(p)$ and $q_a(p)$, where $p=a^2+4b^2$ is a prime such that $p\equiv 1\;\text{mod}\;4$. His congruence reads:
\begin{equation*}\begin{split}
\frac{2E_{p-1}-E_{2p-2}}{p^2}&+\frac{1}{p}\Big(\frac{2E_{p-1}}{p}+4(q_2(p)+q_a(p))+\frac{1}{a^2}\Big)\\&-\frac{1}{a^2}-2(q_2(p)^2+q_a(p)^2)+\frac{3}{8a^4}=0\mpu
\end{split}\end{equation*}

\section{Congruences related to Eulerian numbers}

The recursion formula $(R_{n,m})$ for the Eulerian numbers introduced in $\S\,1$ applied with $n=p-1$ reads:
\begin{equation}E(p-1,m)=(p-1-m)E(p-2,m-1)+(m+1)E(p-2,m)\end{equation}
We start from
\begin{equation}\sum_{m=0}^{p-2}(-1)^m\,E(p-1,m)=0\end{equation}
and use this recursion formula. \\
It yields:
\begin{equation}\begin{split}
0=&\sum_{m=1}^{p-3}(-1)^m\Big(E(p-2,m)-E(p-2,m-1)\Big)\\
&+\sum_{m=1}^{p-3}(-1)^m\Big(E(p-2,m)-E(p-2,m-1)\Big)
\end{split}\end{equation}
After expanding the sums, we obtain:
\begin{equation}\sum_{m=1}^{p-4}(-1)^m(2m+3)E(p-2,m)+p=0\end{equation}
For each integer $k$ with $0\leq k\leq p-3$, we apply the closed formula of $\S\,1.1$ with $n=p-2$. We get:
\begin{equation}
E(p-2,k)=\sum_{s=0}^{k}(-1)^s\binom{p-1}{s}(k+1-s)^{p-2}
\end{equation}
From these identities, we derive the congruences stated in the following lemma.
\newtheorem{Lemma}{Lemma}
\begin{Lemma}
\begin{equation}
\forall 0\leq k\leq p-3,\;E(p-2,k)=H_{k+1}\mpu
\end{equation}
\end{Lemma}
The proof of Lemma $1$ relies on Identity $(5)$ in which we use Fermat's little theorem and the fact that for this whole range of $s$, we have (see Lemma $5$ of \cite{LEVA}):
\begin{equation}
\binom{p-1}{s}=(-1)^s\mpu
\end{equation}
Then, modulo $p$, Equation $(4)$ reduces to:
$$2\sum_{m=1}^{p-4}(-1)^m\,m\,H_{m+1}+3\sum_{m=1}^{p-4}(-1)^m\,H_{m+1}=0\mpu$$
Next, we do the change of indices $k=m+1$ and obtain:
\begin{equation}
-2\sum_{k=2}^{p-3}(-1)^k\,kH_k=\sum_{k=2}^{p-3}(-1)^kH_k\mpu
\end{equation}
We state this intermediate result as a proposition using a more common range for $k$, namely $0\leq k\leq p-1$. For that, we use Wolstenholme's theorem \cite{WOLS} which asserts in particular that $H_{p-1}=0\mpd$. All the other corrective terms cancel each other modulo $p$. Following a usual convention, we set $H_0:=0$.
\newtheorem{Proposition}{Proposition}
\begin{Proposition}
\begin{equation}
-2\sum_{k=0}^{p-1}(-1)^k\,kH_k=\sum_{k=0}^{p-1}(-1)^kH_k\mpu
\end{equation}
\end{Proposition}
Further, by Proposition $5$ of \cite{LEVA}, we also have,
\begin{equation}
\sum_{k=0}^{p-1}H_k=1\,\mpu
\end{equation}
The latter congruence arises from summing all the Eulerian numbers $E(p-2,k)$ for $k$ varying between $0$ and $p-3$, yielding the total number of permutations of $Sym(p-2)$, that is $(p-2)!$. It is then simply a matter of using Lemma $1$ and Wilson's theorem which asserts that $(p-1)!=-1\mpu$, and adding $H_{p-1}$ at the end, which by Wolstenholme's theorem does not contribute \cite{WOLS}. \\
Proposition $1$ and Congruence $(10)$ imply:
\begin{equation}
1-2\sum_{k=0}^{p-1}(-1)^k\,kH_k=2\sum_{\begin{array}{l}k=0\\\text{$k$ even}\end{array}}^{p-1}H_k\mpu
\end{equation}
In order to conclude, we will need a lemma.
\begin{Lemma}
Let $H^{'}_{p-1}$ denote the sum of order $p-1$ of the odd reciprocals of integers.
$$\hp:=1+\frac{1}{3}+\frac{1}{5}+\dots+\frac{1}{p-2}$$
The following congruence holds:
\begin{equation}
\se=\frac{1}{2}(1-\hp)\mpu
\end{equation}
\end{Lemma}
\textsc{Proof of Lemma $2$.}
The proof relies on a result of \cite{LEVA} providing the residue of the Fermat quotient $q_2(p)$ in terms of the residue of twice the number of permutations of $Sym(p-2)$ with an even number of ascents. A reformulation of Theorem $1$ of \cite{LEVA} is:
\begin{equation}
q_2(p)=2N_{p-2}-1\;\mpu,
\end{equation}
where $N_k$ denotes the number of permutations of $Sym(k)$ with an even number of ascents.
The proof of the latter fact is based on several congruences by Emma Lehmer in \cite{LEHM}, on the formula relating an alternating sum of Eulerian numbers to the Bernoulli numbers and on a special case of the Lerch formula \cite{LERC} p. $474$ combined with Wolstenholme's theorem \cite{WOLS}, the latter two providing together the congruence:
\begin{equation}H_{\frac{p-1}{2}}=-2q_2(p)\mpu\end{equation}
Going back to the proof of Lemma $2$, we have:
\begin{eqnarray}
\se&=&1-\so\qquad\qquad\negthickspace\mpu\\
&=&1-\sum_{\begin{array}{l}k=0\\\text{$k$ even}\end{array}}^{p-3}E(p-2,k)\mpu\\
&=&1-N_{p-2}\qquad\qquad\qquad\;\,\mpu\\
&=&\frac{1-q_2(p)}{2}\qquad\qquad\,\,\qquad\mpu
\end{eqnarray}
Further, by Eisenstein's formula for the Fermat quotient (see \cite{EISE}),
\begin{equation}q_2(p)=\frac{1}{2}\sum_{l=1}^{p-1}\frac{(-1)^{l-1}}{l}\mpu\end{equation}
Then, we have by using also Wolstenholme's theorem \cite{WOLS}:
\begin{equation}
q_2(p)=\hp\mpu
\end{equation}
The result stated in Lemma $2$ follows from plugging $(20)$ into $(18)$.\\

The conjunction of Congruence $(11)$ and of Lemma $2$ now yields the following congruence.
\newtheorem{Theorem}{Theorem}
\begin{Theorem}
\begin{equation}
\sum_{k=0}^{p-1}(-1)^k\,kH_k=\frac{1}{2}\,H^{'}_{p-1}\;\mpu
\end{equation}
\end{Theorem}
We derive an immediate corollary.
\newtheorem{Corollary}{Corollary}
\begin{Corollary}
\begin{equation}
\sum_{\begin{array}{l}k=0\\\text{$k$ even}\end{array}}^{p-1}kH_k=\frac{1}{4}(H^{'}_{p-1}-1)\mpu
\end{equation}
\end{Corollary}
\textsc{Proof of Corollary $1$.} By Remark $1.1$ of \cite{SUNB} or Congruence $(9)$ of Corollary $1.3$ of \cite{MESA}, we have:
\begin{equation}
\sum_{k=0}^{p-1}kH_k=-\frac{1}{2}\mpu
\end{equation}

From Identity $(14)$ of \cite{SPIV} which reads:
\begin{equation}
\sum_{k=0}^n\binom{n}{k}H_k=2^n\Big(H_n-\sum_{k=1}^n\frac{1}{k2^k}\Big),
\end{equation}
we derive yet another congruence. Indeed, by setting $n=p-1$ in Spivey's identity, we obtain immediately:
\begin{equation}
\sum_{k=1}^{p-1}\frac{1}{k2^k}=-\sum_{k=0}^{p-1}(-1)^kH_k\;\mpu
\end{equation}
From there, Proposition $1$ and Theorem $1$ imply the following statement.
\begin{Theorem}
\begin{equation}
\sum_{k=1}^{p-1}\frac{1}{k2^k}=\hp\mpu
\end{equation}
\end{Theorem}

In \cite{SPIV}, Michael Spivey offers a new method for evaluating binomial sums in a general setting. He considers sequences $(b_k)_{k\geq 0}$ and $(a_k)_{k\geq 0}$ such that $a_k:=b_{k+1}-b_k$ for $k\geq 0$ and their respective "binomial transforms":
\begin{equation*} h_n=\sum_{k=0}^n\binom{n}{k}b_k\qquad g_n=\sum_{k=0}^n\binom{n}{k}a_k\end{equation*}
As part of his results of \cite{SPIV}, Spivey shows that, where $\delta$ denotes the Kronecker symbol:
\begin{eqnarray}
\forall\,n,\;g_n&=&h_{n+1}-2h_n-b_0\delta_{\lbrace n=-1\rbrace}\\
\forall\,n\geq 0,\;h_n&=&2^n\Big(b_0+\sum_{k=1}^n\frac{g_{k-1}}{2^k}\Big)
\end{eqnarray}
Identity $(27)$ and $(28)$ are respectively Theorem $2$ and $4$ of \cite{SPIV}. In Identity $(27)$, it is understood that the binomial transforms are zero when the integer is negative. \\
Identity $(24)$ is simply obtained by using $(28)$ with $b_k=H_k$ and $a_k=\frac{1}{k+1}$.
\newtheorem{Remark}{Remark}
\begin{Remark}
Proposition $1$ can be derived using $(27)$ with $b_k=kH_k$ and $a_k=1+H_k$. The reduction of
$g_{p-1}=h_p-2h_{p-1}$ modulo $p$ namely provides the result.
\end{Remark}
\begin{Remark}
Proposition $1$ can also be derived yet independently from Peter Paule and Carsten Schneider's identity
\begin{equation*}
\sum_{j=0}^n(n-2j)H_j\binom{n}{j}=1-2^n\qquad n\geq 0
\end{equation*}
reduced modulo $p$ with $n=p-1$. Their result from \cite{PAUL} is based in particular on Zeilberger's algorithm for definite hypergeometric sums, see \cite{ZEIL}, and on expressing harmonic numbers in terms of differentiation of binomial coefficients, a method which can be traced back to Isaac Newton.
Some of these types of sums, more widely studied by the authors, closely relate to other binomial sums playing an important role in Ap\'ery's approach to prove the irrationality of $\zeta(2)$ and $\zeta(3)$, see \cite{APER}.
\end{Remark}
\begin{Remark}
Congruence $(26)$ was first uncovered by J.W.L. Glaisher in \cite{GLAI} under the form:
\begin{equation*}
\sum_{k=1}^{p-1}\frac{2^k}{k}=-2q_2(p)\mpu
\end{equation*}
The new proof presented here reflects a different viewpoint centered on the Eulerian numbers.
In \cite{SUNH}, using this time integration methods, Zhi-Hong Sun also shows that:
\begin{equation*}
\sum_{k=1}^{p-1}\frac{2^k}{k^2}=-q_2(p)^2+p\Big(\frac{2}{3}q_2(p)^3+\frac{7}{6}B_{p-3}\Big)\;\mpt,
\end{equation*}
thus generalizing Andrew Granville's congruence modulo $p$, see \cite{GRAN}.
\end{Remark}
\begin{Remark}
Congruences like in Theorem $1$ or Corollary $1$ involving more generally powers of integers equal to or greater than one and generalized harmonic numbers or their powers is the purpose of ongoing research, see in chronological order \cite{MESA}, \cite{WAY1} and \cite{WAY2}. Only the non alternating sums have been investigated so far.
\end{Remark}
\section{Congruences related to Euler numbers}
In what follows, $\mg_k$ denotes the $k$-th divided Genocchi number. \\
In this part, we will need two lemmas, both involving the Genocchi numbers. These two lemmas get both derived from Kwang-Wu Chen's closed formula respectively applied with $n=p-1$ and $n=p-2$. We start with $n=p-1$. \\
We have:
\begin{equation}
E_{p-1}=1+\frac{1}{p}\sum_{k=2}^{p}\binom{p}{k}2^{k-1}G_k
\end{equation}
Moreover,
$$\forall\,2\leq k\leq p-1,\;\binom{p}{k}=\binom{p-1}{k}\frac{p}{p-k}\equiv\frac{(-1)^{k+1}}{k}\,p\;\mpd$$
By using the original congruence by George Stetson Ely, we thus get:
\begin{Lemma} Let the $\mg_k$'s denote the divided Genocchi numbers.
\begin{equation}
\sum_{k=2}^{p-1}2^k\mg_k=\begin{cases}2\;\;\;\;\mpu&\text{if $p\equiv 1\;\text{mod}\;4$}\\
-2\;\mpu&\text{if $p\equiv 3\;\text{mod}\;4$}\end{cases}
\end{equation}
\end{Lemma}
\begin{Corollary} Let the $\mg_k$'s denote the divided Genocchi numbers.
\begin{equation}
\sum_{k=2}^{p-3}2^k\mg_k=\begin{cases}4\,\left(\begin{array}{l}\#\text{permutations on $(p-2)$ letters}\\\text{with an even number of ascents}\end{array}\right)\;\mpu&\text{if $p\equiv 1\;\text{mod}\;4$}\\
4\,\left(\begin{array}{l}\#\text{permutations on $(p-2)$ letters}\\\text{with an even number of ascents}\\\text{and distinct from the identity}\end{array}\right)\;\mpu&\text{if $p\equiv 3\;\text{mod}\;4$}\end{cases}
\end{equation}
\end{Corollary}
\textsc{Proof of the corollary.} We have:
\begin{equation}G_{p-1}=2(1-2^{p-1})B_{p-1}=-2q_2(p)pB_{p-1}\equiv 2q_2(p)\;\mpu\end{equation}
The last congruence holds as by the Von Staudt-Clausen's theorem \cite{VONS}\cite{CLAU}, we have $pB_{p-1}=-1\mpu$. It then suffices to apply Lemma $3$ above and Theorem $1$ of \cite{LEVA} in order to conclude. \\

We continue with $n=p-2$ in the closed formula. It yields:
\begin{equation}
E_{p-2}=1+\frac{1}{p-1}\sum_{k=2}^{p-1}\binom{p-1}{k}2^{k-1}G_k
\end{equation}
We obtain:
\begin{Lemma}
\begin{equation}
\sum_{k=2}^{p-1}2^kG_k=2\;\mpu
\end{equation}
\end{Lemma}
Because the Genocchi numbers with odd subscripts are zero, the congruence of Lemma $4$ is equivalent to the one below.
\begin{equation}
\sum_{\begin{array}{l}1\leq k\leq p-2\\\;\;\;\;\;\text{k odd}\end{array}}2^kG_{k+1}=1\mpu
\end{equation}
Moreover, when $k$ is odd, we have:
\begin{equation}
2^kG_{k+1}=-(k+1)\widehat{E_k}
\end{equation}
Thus,
\begin{equation}
\sum_{\begin{array}{l}1\leq k\leq p-2\\\;\;\;\;\;\;\text{k odd}\end{array}}(k+1)\wh{E_k}=-1\mpu
\end{equation}
Further, starting again from $(36)$, we notice that:
\begin{eqnarray}
\sum_{\begin{array}{l}1\leq k\leq p-2\\\;\;\;\;\;\;\text{k odd}\end{array}}\wh{E_k}&=&-\sum_{\begin{array}{l}1\leq k\leq p-2\\\;\;\;\;\;\;\text{k odd}\end{array}}2^k\mg_{k+1}\\
&=&-\frac{1}{2}\sum_{k=2}^{p-1}2^k\mg_k
\end{eqnarray}
Then,
\begin{equation}
\sum_{\begin{array}{l}1\leq k\leq p-2\\\;\;\;\;\;\;\text{$k$ odd}\end{array}}k\wh{E_k}=-1+\frac{1}{2}\sum_{k=2}^{p-1}2^k\mg_k\;\mpu
\end{equation}
Furthermore,
\begin{equation}
T_{2k+1}=(-1)^k\wih{E_{2k+1}}
\end{equation}
By applying Lemma $3$, we thus derive the following theorem.
\begin{Theorem} Let the $T_m$'s denote the tangent numbers. The following congruence holds.
\begin{equation*}
\sum_{\begin{array}{l}m=1\\\text{$m$ odd}\end{array}}^{p-1}(-1)^{\frac{m+1}{2}}mT_m=\begin{cases}0\;\mpu&\text{if $p\equiv\,1$ mod 4}\\
2\;\mpu&\text{if $p\equiv\,3$ mod 4}\end{cases}
\end{equation*}
\end{Theorem}
Note that a reformulation of Theorem $3$ is the following.
\begin{equation*}
\sum_{k=1}^{\frac{p-1}{2}}(2k-1)\widehat{E_{2k-1}}=\begin{cases}0\;\;\;\;\;\text{mod}\,p&\text{if $p\equiv 1\;\text{mod}\;4$}\\
-2\;\;\text{mod}\,p&\text{if $p\equiv 3\;\text{mod}\;4$}\end{cases}\end{equation*}
Lemma $4$ could also be derived using the Yuan He and Qunying Liao congruence combined with the Lehmer-Mirimanoff congruences which we recall below.
\begin{Theorem} Due to Emma Lehmer in $1938$ for $k$ even \cite{LEHM} and due to Dmitry Mirimanoff in $1895$ \cite{MIRI} for $k$ odd.
Assume that $p-1\not|k-1$. Then,
\begin{equation}
S_k:=\sum_{r=1}^{\frac{p-1}{2}}r^k=\begin{cases}\big(2^{-k+1}-1\big)\frac{p}{2}B_k\;\;\mpd &\text{if $k$ is even}\\&\\\big(\frac{1}{2^{k+1}}-1\big)\frac{2B_{k+1}}{k+1}\;\;\mpd&\text{if $k$ is odd}
\end{cases}
\end{equation}
\end{Theorem}
Applied with $m=p$ and $n=p-2$, Y. He and Q.Y. Liao's congruence reads:
\begin{equation}
\sum_{j=0}^{p-1}(-1)^j(2j+1)^{p-2}=0\mpu
\end{equation}
If we expand using the binomial formula, then use some classical congruences on the binomial coefficients, some congruences on the sums of powers of the first $p-1$ integers and the Lehmer--Mirimanoff congruences, then it leads to the result stated in Lemma $4$. \\
Interestingly, the old Mirimanoff congruences imply another congruence similar to those of Lemmas $3$ and $4$, involving powers of two and Genocchi numbers. \\Indeed, by summing the Mirimanoff congruences over the range of $k$ corresponding to $3\leq k\leq p-2$, we obtain:
\begin{equation}
\sum_{k=4}^{p-1}\frac{\mg_k}{2^k}=S_3+S_5+\dots+S_{p-2}\;\mpd
\end{equation}
It is shown in \cite{LEVA} that the sum of the first $\frac{p-1}{2}$ odd powers of the first $\frac{p-1}{2}$ integers is congruent to $-\frac{1}{2}$ modulo $p$. This is Theorem $7$ of \cite{LEVA}. Expressed in terms of the $S_k$'s, we thus have:
\begin{equation}
S_1+S_3+S_5+\dots+S_{p-2}=-\frac{1}{2}\mpu
\end{equation}
We deduce:
\begin{Theorem} Let $G_k$ denote the $k$-th Genocchi number. The following congruence holds.
\begin{equation*}
\sum_{k=1}^{p-1}\frac{G_k}{k2^k}=0\mpu
\end{equation*}
\end{Theorem}
We note that similar congruences had been worked out for harmonic numbers and Bernoulli numbers. In \cite{SUNB}, Zhi-Wei Sun showed that
\begin{equation}
\sum_{k=1}^{p-1}\frac{H_k}{k2^k}=0\mpu
\end{equation}
In \cite{LEVA}, letting $\mb_k:=\frac{B_k}{k}$ when $k>0$ and $\mb_0:=-\frac{1}{p}$, it is shown that:
\begin{equation}
\sum_{k=0}^{p-1}\frac{\mb_k}{2^k}=0\;\mpu,
\end{equation}
or equivalently, excluding the upper and lower bounds of the sum:
\begin{equation}
\sum_{k=1}^{p-2}\frac{B_k}{k2^k}=-\frac{1}{2}H_{\frac{p-1}{2}}+\frac{pB_{p-1}+1}{p}-1\mpu
\end{equation}
These are respectively Theorem $5$ and Congruence $(11)$ in Theorem $6$ of \cite{LEVA}. \\
More congruences of similar type involving harmonic numbers or generalized harmonic numbers appear in \cite{SUNZ} and in Theorem $1.1$ of \cite{MESS}.\\\\

\textbf{Acknowledgements.} I thank Peter Moree for drawing my attention on the Genocchi numbers. \\

\textsc{Email address:} \textit{clairelevaillant@yahoo.fr}

\end{document}